\documentclass{amsart}
\usepackage{amsfonts,amsmath,amsthm,amssymb}
\newtheorem{theorem}{Theorem}
\newtheorem{corollary}{Corollary}
\newtheorem{proposition}{Proposition}
\title[Integral group ring of the McLaughlin simple group]
{Integral group ring of the \\ McLaughlin simple group}
\date{}

\author{V.A.~Bovdi, A.B.~Konovalov}

\address{V.A.~Bovdi
\newline Institute of Mathematics, University of Debrecen,
\newline P.O.  Box 12, H-4010 Debrecen, Hungary
\newline Institute of Mathematics and Informatics, College of Ny\'\i regyh\'aza,
\newline S\'ost\'oi \'ut 31/b, H-4410 Ny\'\i regyh\'aza, Hungary}
\email{vbovdi@math.klte.hu}

\address{A.B.~Konovalov
\newline School of Computer Science, University of St Andrews,
\newline Jack Cole Building, North Haugh, St Andrews, Fife, KY16 9SX, Scotland}
\email{konovalov@member.ams.org}

\thanks{The research was supported by OTKA grants No.T 43034, No.K61007}

\dedicatory{Dedicated to 60th birthday of Professor V.I.~Sushchansky}


\thanks{}


\begin{document}
\begin{abstract}
We consider the Zassenhaus conjecture for the normalized unit
group of the integral group ring of the McLaughlin sporadic group
$\verb"McL"$. As a consequence, we confirm for this group the 
Kimmerle's conjecture on prime graphs.
\end{abstract}

\maketitle

\section{Introduction, conjectures and main results}
\label{Intro}

Let $V(\mathbb Z G)$ be  the normalized unit group of the
integral group ring $\mathbb Z G$ of  a finite group $G$. A long-standing
conjecture  of H.~Zassenhaus {\bf (ZC)} says that every torsion unit
$u\in V(\mathbb ZG)$ is conjugate within the rational group algebra
$\mathbb Q G$ to an element in $G$ (see \cite{Zassenhaus}).

For finite simple groups the main tool for the investigation of the
Zassenhaus conjecture is the Luthar--Passi method, introduced in
\cite{Luthar-Passi} to solve it for $A_{5}$. 
Later M.~Hertweck in \cite{Hertweck1} extended
the Luthar--Passi method and applied it for the investigation
of the Zassenhaus conjecture for $PSL(2,p^{n})$.
The Luthar--Passi method proved to be useful for groups containing 
non-trivial normal subgroups as well. For some recent results we refer to
\cite{Bovdi-Hofert-Kimmerle,Bovdi-Konovalov, Hertweck2,
Hertweck1, Hertweck3, Hofert-Kimmerle}. 
Also, some related properties and  some weakened variations of the
Zassenhaus  conjecture can be found in
\cite{Artamonov-Bovdi,Luthar-Trama} and
\cite{Bleher-Kimmerle,Kimmerle}.

First of all, we need to introduce some notation. By $\# (G)$ we
denote the set of all primes dividing the order of $G$. The
Gruenberg--Kegel graph (or the prime graph) of $G$ is the graph
$\pi (G)$ with vertices labeled by the primes in $\# (G)$ and with
an edge from $p$ to $q$ if there is an element of order $pq$ in
the group $G$. In \cite{Kimmerle} W.~Kimmerle   proposed the
following weakened variation of the Zassenhaus conjecture:

\begin{itemize}
\item[]{\bf (KC)} \qquad
If $G$ is a finite group then $\pi (G) =\pi (V(\mathbb Z G))$.
\end{itemize}

In particular, in the same  paper  W.~Kimmerle verified   that
{\bf (KC)} holds for finite Frobenius and solvable groups. We
remark that with respect to the so-called $p$-version of the Zassenhaus 
conjecture the investigation
of Frobenius groups was completed by  M.~Hertweck and the first
author  in \cite{Bovdi-Hertweck}. In \cite{Bovdi-Jespers-Konovalov,
Bovdi-Konovalov, Bovdi-Konovalov-M24, Bovdi-Konovalov-M23, Bovdi-Konovalov-Linton, 
Bovdi-Konovalov-Siciliano} 
{\bf (KC)} was confirmed for the Mathieu simple groups $M_{11}$, 
$M_{12}$, $M_{22}$, $M_{23}$, $M_{24}$ and the sporadic Janko simple groups 
$J_1$, $J_{2}$ and $J_3$.

Here  we continue these investigations for the McLaughlin simple group
$\verb"McL"$. Although using the Luthar--Passi method we cannot prove the
rational conjugacy for torsion units of $V(\mathbb Z \verb"McL")$, our
main result gives a lot of information on partial augmentations of
these units. In particular, we confirm the Kimmerle's conjecture for
this group.

Let $G=\verb"McL"$. It is well known
(see \cite{GAP}) that 
$|G|=2^7 \cdot 3^6 \cdot 5^3 \cdot 7 \cdot 11$ and
$exp(G)=2^3 \cdot 3^2 \cdot  5 \cdot 7 \cdot 11$.
Let
\[
\begin{split}
\mathcal{C} =\{
C_{1},C_{2a},C_{3a}, & C_{3b},C_{4a},C_{5a},C_{5b},C_{6a},C_{6b},
  C_{7a},C_{7b},C_{8a},C_{9a},C_{9b},\\
  & C_{10a}, C_{11a},C_{11b},C_{12a},C_{14a},C_{14b},C_{15a},C_{15b},C_{30a},C_{30b} \}
\end{split}
\]
be the collection of all conjugacy classes of $\verb"McL"$, where the first
index denotes the order of the elements of this conjugacy class
and $C_{1}=\{ 1\}$. Suppose $u=\sum \alpha_g g \in V(\mathbb Z G)$
has finite order $k$. Denote by
$\nu_{nt}=\nu_{nt}(u)=\varepsilon_{C_{nt}}(u)=\sum_{g\in C_{nt}}
\alpha_{g}$ the partial augmentation of $u$ with respect to
$C_{nt}$. From the Berman--Higman Theorem 
(see \cite{Berman} and \cite{Sandling}, Ch.5, p.102)
one knows that
$\nu_1 =\alpha_{1}=0$ and
\begin{equation}\label{E:1}
\sum_{C_{nt}\in \mathcal{C}} \nu_{nt}=1.
\end{equation}
Hence, for any character $\chi$ of $G$, we get that $\chi(u)=\sum
\nu_{nt}\chi(h_{nt})$, where $h_{nt}$ is a representative of the
conjugacy class $ C_{nt}$.

Our main result is the following

\begin{theorem}\label{T:1}

Let $G$ denote the McLaughlin simple group $\verb"McL"$. Let  $u$ be
a torsion unit of $V(\mathbb ZG)$ of order $|u|$. Denote by
$\frak{P}(u)$ the tuple
\[
\begin{split}
( \nu_{2a},\nu_{3a},\nu_{3b},&\nu_{4a},\nu_{5a},\nu_{5b},\nu_{6a},\nu_{6b},
  \nu_{7a},\nu_{7b},\nu_{8a},\nu_{9a},\nu_{9b},\\
  & \nu_{10a},\nu_{11a},\nu_{11b},\nu_{12a},\nu_{14a},\nu_{14b},\nu_{15a},
    \nu_{15b},\nu_{30a},\nu_{30b})\in \mathbb Z^{23}
\end{split}
\]
of partial augmentations of $u$ in $V(\mathbb ZG)$. The following
properties hold.

\begin{itemize}

\item[(i)] There is no elements of orders $21$, $22$, $33$, $35$, $55$, $77$ in 
$V(\mathbb ZG)$. Equivalently, if $|u|\not\in\{18,20,24,28,36,40,45,56,60,72,90,120,180,360\}$, then $|u|$ 
coincides with the order of some element $g\in G$.

\item[(ii)] If $|u|=2$, then $u$ is rationally conjugate to some $g\in G$.

\item[(iii)]
If $|u|=3$, then all components of $\frak{P}(u)$ are zero except possibly 
$\nu_{3a}$ and $\nu_{3b}$, and the pair $(\nu_{3a},\nu_{3b})$ is one of 
$$
 \{ \;  ( -2, 3 ), \; ( -1, 2 ), \; ( 0, 1 ), \; ( 1, 0 ) \; \} .
$$

\item[(iv)]
If $|u|=5$, then all components of $\frak{P}(u)$ are zero except possibly 
$\nu_{5a}$ and $\nu_{5b}$, and the pair $(\nu_{5a},\nu_{5b})$ is one of 
\[
\begin{split}
 \{ \; ( -4, 5 ), \; ( -3, 4 ), \; ( -2, 3 ), \; ( -1, 2 ), \; ( 0, 1 ), \; ( 1, 0 ) \; \} .
\end{split}
\]

\item[(v)]
If $|u|=7$, then all components of $\frak{P}(u)$ are zero except possibly 
$\nu_{7a}$ and $\nu_{7b}$, and the pair $(\nu_{7a},\nu_{7b})$ is one of 
$$
 \{ \; (\nu_{7a},\nu_{7b}) \; \mid \; -86 \leq \nu_{7a} \leq 87, \quad \nu_{7a}+\nu_{7b}=1 \; \} .
$$

\item[(vi)]
If $|u|=11$, then all components of $\frak{P}(u)$ are zero except possibly 
$\nu_{11a}$ and $\nu_{11b}$, and the pair $(\nu_{11a},\nu_{11b})$ is one of 
$$
 \{ \; (\nu_{11a},\nu_{11b}) \; \mid \; -9 \leq \nu_{11a} \leq 10, \quad \nu_{11a}+\nu_{11b}=1 \; \} .
$$

\end{itemize}
\end{theorem}

As an immediate consequence of  part (i) of the Theorem we obtain

\begin{corollary} If $G=\verb"McL"$ then
$\pi(G)=\pi(V(\mathbb ZG))$.
\end{corollary}

\section{Preliminaries}
The following result is a reformulation of
the Zassenhaus conjecture in terms of vanishing of
partial augmentations of torsion units.

\begin{proposition}\label{P:5}
(see \cite{Luthar-Passi} and
Theorem 2.5 in \cite{Marciniak-Ritter-Sehgal-Weiss})
Let $u\in V(\mathbb Z G)$
be of order $k$. Then $u$ is conjugate in $\mathbb
QG$ to an element $g \in G$ if and only if for
each $d$ dividing $k$ there is precisely one
conjugacy class $C$ with partial augmentation
$\varepsilon_{C}(u^d) \neq 0 $.
\end{proposition}

The next result now yield that several partial augmentations
are zero.

\begin{proposition}\label{P:4}
(see \cite{Hertweck2}, Proposition 3.1;
\cite{Hertweck1}, Proposition 2.2)
Let $G$ be a finite
group and let $u$ be a torsion unit in $V(\mathbb
ZG)$. If $x$ is an element of $G$ whose $p$-part,
for some prime $p$, has order strictly greater
than the order of the $p$-part of $u$, then
$\varepsilon_x(u)=0$.
\end{proposition}

The key restriction on partial augmentations is given 
by the following result that is the cornerstone of
the Luthar--Passi method.

\begin{proposition}\label{P:1}
(see \cite{Luthar-Passi,Hertweck1}) Let either $p=0$ or $p$ a prime
divisor of $|G|$. Suppose
that $u\in V( \mathbb Z G) $ has finite order $k$ and assume $k$ and
$p$ are coprime in case $p\neq 0$. If $z$ is a complex primitive $k$-th root
of unity and $\chi$ is either a classical character or a $p$-Brauer
character of $G$, then for every integer $l$ the number
\begin{equation}\label{E:2}
\mu_l(u,\chi, p ) =
\textstyle\frac{1}{k} \sum_{d|k}Tr_{ \mathbb Q (z^d)/ \mathbb Q }
\{\chi(u^d)z^{-dl}\}
\end{equation}
is a non-negative integer.
\end{proposition}

Note that if $p=0$, we will use the notation $\mu_l(u,\chi
, * )$ for $\mu_l(u,\chi , 0)$.

Finally, we shall use the well-known bound for
orders of torsion units.

\begin{proposition}\label{P:2}  (see  \cite{Cohn-Livingstone})
The order of a torsion element $u\in V(\mathbb ZG)$
is a divisor of the exponent of $G$.
\end{proposition}

\section{Proof of the Theorem}

Throughout this section we denote $\verb"McL"$ by $G$. 
The character table of $G$,
as well as the $p$-Brauer character tables, which will be denoted by
$\mathfrak{BCT}{(p)}$ where $p\in\{2,3,5,7,11\}$, can be found using
the computational algebra system GAP \cite{GAP},
which derives these data from \cite{AFG,ABC}. For the characters
and conjugacy classes we will use throughout the paper the same
notation, indexation inclusive, as used in the GAP Character Table Library.

First of all we will investigate units of orders $2$, $3$, $5$, $7$ and
$11$, since the group $G$ possesses elements of these orders. After this, 
by Proposition \ref{P:2}, the order of each torsion unit divides the 
exponent of $G$, so to prove the Kimmerle's conjecture, it remains to 
consider units of orders $21$, $22$, $33$, $35$, $55$ and $77$. We prove 
that no units of all these orders do  appear in $V(\mathbb ZG)$.

Now we consider each case separately.

\noindent$\bullet$ Let $u$ be an involution. Using
Proposition \ref{P:4} we obtain that all partial
augmentations except one are zero. Thus by Proposition \ref{P:5} the
proof of part (ii) of  Theorem \ref{T:1} is done.

\noindent $\bullet$ Let $u$ be a unit of order $3$. By (\ref{E:1})
and Proposition \ref{P:4} we get $\nu_{3a}+\nu_{3b}=1$.
Put $t_1 = 5 \nu_{3a} - 4 \nu_{3b}$. 
By (\ref{E:2}) we obtain  the system of inequalities
\[
\begin{split} 
\mu_{0}(u,\chi_{2},*) & = \textstyle \frac{1}{3} (-2 t_1 + 22) \geq 0; \qquad 
\mu_{1}(u,\chi_{2},*)   = \textstyle \frac{1}{3} (   t_1 + 22) \geq 0, \\ 
\end{split} 
\]
from which $t_1 \in \{ -22, -19, -16, -13, -10, -7, -4, -1, 2, 5, 8, 11 \}$.
Now for each possible value of $t_1$ consider the system of linear equations
$$ 
\nu_{3a}+\nu_{3b}=1, \qquad 
5 \nu_{3a} - 4 \nu_{3b} = t_1.
$$
Since $ \tiny{ \left|
\begin{matrix}
1&1\\
5&-4\\
\end{matrix}\right| } \not=0$, this system always has the unique 
solution. First we select only integer solutions, and then 
using the condition that all $\mu_i(u,\chi_{j},*)$ are non-negative 
integers, we obtain only four pairs $(\nu_{3a},\nu_{3b})$ listed in 
part (iii) of Theorem \ref{T:1}.

\noindent $\bullet$ Let $u$ be a unit of order $5$. By (\ref{E:1})
and Proposition \ref{P:4} we get $\nu_{5a}+\nu_{5b}=1$. 
Put $t_1 = 3 \nu_{5a} - 2 \nu_{5b}$.
By
(\ref{E:2}) we obtain the system of inequalities
\[
\begin{split} 
\mu_{0}(u,\chi_{2},*) & = \textstyle \frac{1}{5} (-4 t_1 + 22) \geq 0; \qquad 
\mu_{1}(u,\chi_{2},*)   = \textstyle \frac{1}{5} (   t_1 + 22) \geq 0, \\ 
\end{split} 
\]
so $t_1 \in \{ -22, -17, -12, -7, -2, 3 \}$.
Using the same arguments as in the previous case, 
we obtain only six pairs $(\nu_{5a},\nu_{5b})$ 
listed in part (iv) of Theorem \ref{T:1}.

\noindent $\bullet$ Let $u$ be a unit of order $7$. By (\ref{E:1})
and Proposition \ref{P:4} we get $\nu_{7a}+\nu_{7b}=1$. 
Put $t_1 = 4 \nu_{7a} - 3 \nu_{7b}$.
Using $\mathfrak{BCT}{(3)}$ and $\mathfrak{BCT}{(5)}$,
by (\ref{E:2}) we have 
\[
\begin{split} 
\mu_{3}(u,\chi_{7},3) & = \textstyle \frac{1}{7} (t_1 + 605) \geq 0; \quad 
\mu_{1}(u,\chi_{12},5)  = \textstyle \frac{1}{7} (-t_1 + 3245) \geq 0; \\ 
& \mu_{1}(u,\chi_{7},3) = \textstyle \frac{1}{7} (-3 \nu_{7a} + 4 \nu_{7b} + 605) \geq 0. \\ 
\end{split} 
\]
It follows that we have only 174 pairs $(\nu_{7a},\nu_{7b})$, 
given in part (v) of Theorem~\ref{T:1}. Note that using our 
implementation of the Luthar--Passi method, which we intended
to make available in the GAP package LAGUNA \cite{LAGUNA},
we checked that it is not possible to further reduce the number of solutions, and the same
remark also applies for the remaining part of the paper.

\noindent $\bullet$ Let $u$ be a unit of order $11$. By (\ref{E:1})
and Proposition \ref{P:4} we have $\nu_{11a}+\nu_{11b}=1$. Using
$\mathfrak{BCT}{(3)}$, by (\ref{E:2}) we obtain the system of inequalities
\[
\begin{split} 
\mu_{1}(u,\chi_{3},3) & = \textstyle \frac{1}{11} (6 \nu_{11a} - 5 \nu_{11b} + 104) \geq 0; \\ 
\mu_{2}(u,\chi_{3},3) & = \textstyle \frac{1}{11} (-5 \nu_{11a} + 6 \nu_{11b} + 104) \geq 0, \\ 
\end{split} 
\]
that has only that twenty pairs $(\nu_{11a},\nu_{11b})$ 
listed in part (vi) of the Theorem \ref{T:1}.

\noindent$\bullet$ Let $u$ be a unit of order $21$. By (\ref{E:1})
and Proposition \ref{P:4} we have
$$
\nu_{3a}+\nu_{3b}+\nu_{7a}+\nu_{7b}=1.
$$
Put $t_1 = 5 \nu_{3a} - 4 \nu_{3b} -  \nu_{7a} -  \nu_{7b}$, $t_2
= 5 \nu_{3a} + 2 \nu_{3b}$ and $t_3 =3 \nu_{7a} - 4 \nu_{7b}$.
Since  $|u^{7}|=3$, for any character $\chi$ of $G$ we need to
consider four cases, defined by part (iii)  of the Theorem.
Now we consider each case separately:

Case 1. Let $\chi(u^{7}) = \chi(3a)$. Using (\ref{E:2}),
we obtain the system of inequalities
\[
\mu_{3}(u,\chi_{2},*) = \textstyle \frac{1}{21} (2t_1 + 11) \geq
0; \qquad \mu_{0}(u,\chi_{2},*) = \textstyle \frac{1}{21} (-12t_1
+ 18) \geq 0,
 \]
which has no integral solution.

Case 2. Let $\chi(u^{7}) = \chi(3b)$. Again, using (\ref{E:2}),
we obtain the system of inequalities
\[
\begin{split}
\mu_{0}(u,\chi_{2},*) & = \textstyle \frac{1}{21} (-12 t_1 + 36)
\geq 0; \qquad \mu_{7}(u,\chi_{2},*)  = \textstyle \frac{1}{21} (6t_1 + 24) \geq 0; \\ 
\mu_{0}(u,\chi_{3},*) & = \textstyle \frac{1}{21} (36  t_2  + 243) \geq 0; \qquad \mu_{7}(u,\chi_{3},*)= \textstyle \frac{1}{21} (-18 t_2  + 225) \geq 0; \\ 
\mu_{1}(u,\chi_{16},*) & = \textstyle \frac{1}{21} (- t_3 + 8386)
\geq 0; \qquad \mu_{9}(u,\chi_{16},*) = \textstyle \frac{1}{21} (2t_3 + 8386) \geq 0; \\ 
&\quad \mu_{1}(u,\chi_{5},*) = \textstyle \frac{1}{21} (-13 \nu_{3a} + 5 \nu_{3b} + 765) \geq 0. \\ 
\end{split}
\]
This yields 
$t_1 \in \{ -4, 3 \}$, 
$t_2 \in \{-5, 2, 9\}$ and 
$t_3 \in \{ 7+21k \mid -200 \leq k \leq 399 \}$,
but none of possible combinations of $t_i$'s gives
us any solution.

Case 3. Let $\chi(u^{7}) = -2\chi(3a)+3\chi(3b)$. Then using
(\ref{E:2}), we obtain the system
\[
\begin{split}
\mu_{1}(u,\chi_{2},*) & = \textstyle \frac{1}{21} (-t_1 - 1) \geq
0; \qquad\quad
\mu_{7}(u,\chi_{2},*)  = \textstyle \frac{1}{21} (6t_1 + 6) \geq 0; \\ 
\mu_{0}(u,\chi_{3},*) & = \textstyle \frac{1}{21} (36t_2 + 207)
\geq 0; \qquad
\mu_{7}(u,\chi_{3},*)  = \textstyle \frac{1}{21} (-18t_2 + 243) \geq 0; \\ 
\mu_{1}(u,\chi_{16},*) & = \textstyle \frac{1}{21} (-t_3 + 8218)
\geq 0; \qquad
\mu_{9}(u,\chi_{16},*)  = \textstyle \frac{1}{21} (2t_3 + 8218) \geq 0, \\ 
\end{split}
\]
from which $t_1=-1$, $t_2 \in \{ -4, 3, 10 \}$ and 
$t_3 \in \{ 7+21k \mid -196 \leq k \leq 391 \}$, 
and again we have no solution for every combination of $t_i$'s.

Case 4. Let $\chi(u^{7}) = -\chi(3a)+2\chi(3b)$.
Using (\ref{E:2}), we obtain the system
\[
\begin{split}
\mu_{7}(u,\chi_{2}, *) & = \textstyle \frac{1}{21} (6t_1 + 15)
\geq 0; \qquad\quad
\mu_{0}(u,\chi_{2},*)  = \textstyle \frac{1}{21} (-12t_1 + 54) \geq 0; \\ 
\mu_{0}(u,\chi_{3},*) & = \textstyle \frac{1}{21} (36t_2 + 225)
\geq 0; \qquad
\mu_{7}(u,\chi_{3},*)  = \textstyle \frac{1}{21} (-18t_2 + 234) \geq 0; \\ 
\mu_{9}(u,\chi_{16},*) & = \textstyle \frac{1}{21} (2t_3 + 8015)
\geq 0; \qquad
\mu_{1}(u,\chi_{16},*) = \textstyle \frac{1}{21} (- t_3 + 8015) \geq 0, \\ 
\end{split}
\]
so $t_1 =1$, $t_2 \in\{  -1, 6, 13 \}$ and
$t_3 \in \{ 14+21k \mid -191 \leq k \leq 381 \}$,
that also gives us no solutions.


\noindent$\bullet$ Let $u$ be a unit of order $22$. 
By (\ref{E:1})
and Proposition \ref{P:4} we have
$$
\nu_{2a}+\nu_{11a}+\nu_{11b}=1.
$$
Now by (\ref{E:2}) we obtain the system of inequalities 
\[
\mu_{0}(u,\chi_{2},*) = \textstyle \frac{1}{22} (60 \nu_{2a} + 28)
\geq 0; \qquad \mu_{11}(u,\chi_{2},*)  = \textstyle \frac{1}{22}
(-60 \nu_{2a} + 16) \geq 0,
\]
which has no integral solution.

\noindent$\bullet$ Let $u$ be a unit of order $33$. By (\ref{E:1})
and Proposition \ref{P:4} we have 
$$
\nu_{3a}+\nu_{3b}+\nu_{11a}+\nu_{11b}=1.
$$
Put $t_1 = 5 \nu_{3a} - 4 \nu_{3b}$, 
$t_2 = 5 \nu_{3a} + 2 \nu_{3b}$ and 
$t_3=32 \nu_{3a} - 4 \nu_{3b} - 6 \nu_{11a} + 5 \nu_{11b}$.
Since  $|u^{11}|=3$, for any character $\chi$ of $G$ we need to
consider four   cases, defined by part (iii) of the Theorem.

Case 1. Let $\chi(u^{11}) = \chi(3a)$. Then
by (\ref{E:2}) we obtain the system of inequalities 
\[
\mu_{11}(u,\chi_{2},*) = \textstyle \frac{1}{33} (10t_1 + 27) \geq
0; \qquad \mu_{0}(u,\chi_{2},*) = \textstyle \frac{1}{33} (-20t_1 + 12) \geq 0, \\ 
\]
that has no integral solution.

Case 2. Let $\chi(u^{11}) = \chi(3b)$. Now (\ref{E:2}) gives us the system
\[
\mu_{11}(u,\chi_{2},*) = \textstyle \frac{1}{33} (10t_1 + 18) \geq
0; \qquad
\mu_{0}(u,\chi_{2},*)  = \textstyle \frac{1}{33} (-20 t_1 + 30) \geq 0, \\ 
\]
which also has no integral solution.

Case 3. Let $\chi(u^{11}) = -2\chi(3a)+3\chi(3b)$.
By (\ref{E:2}) we obtain that
\[
\begin{split}
\mu_{1}(u,\chi_{2},*) & = \textstyle \frac{1}{33} (-5 \nu_{3a} + 4
\nu_{3b} ) \geq 0; \qquad
\mu_{11}(u,\chi_{2},*)  = \textstyle \frac{1}{33} (50 \nu_{3a} - 40 \nu_{3b} ) \geq 0; \\ 
\mu_{0}(u,\chi_{3},*) & = \textstyle \frac{1}{33} (60t_2 + 207)
\geq 0; \qquad \quad
\mu_{11}(u,\chi_{3},*) = \textstyle \frac{1}{33} (-30t_2 + 243) \geq 0; \\ 
\mu_{1}(u,\chi_{7},*) & = \textstyle \frac{1}{33} (t_3 + 978) \geq
0; \qquad \quad \quad
\mu_{3}(u,\chi_{7},*)  = \textstyle \frac{1}{33} (-2t_3 + 750) \geq 0. \\ 
\end{split}
\]
It follows that $t_1=0$, $t_2=7$ and
$t_3 \in \{ 12+33k \mid -30 \leq k \leq 11 \}$, and we have no solutions again.

Case 4. Let $\chi(u^{11}) = -\chi(3a)+2\chi(3b)$.
By (\ref{E:2}) we obtain the system
\[
\mu_{11}(u,\chi_{2},*)  = \textstyle \frac{1}{33} (10t_1 + 9) \geq 0; \qquad 
\mu_{0}(u,\chi_{2},*)  = \textstyle \frac{1}{33} (-20t_1 + 48) \geq 0, \\ 
\]
which has no integral solution.

\noindent$\bullet$ Let $u$ be a unit of order $35$. By (\ref{E:1})
and Proposition \ref{P:4} we have
$$
\nu_{5a}+\nu_{5b}+\nu_{7a}+\nu_{7b}=1.
$$
Put $t_1 = 3 \nu_{5a} - 2 \nu_{5b} -  \nu_{7a} -  \nu_{7b}$, $t_2
= 6 \nu_{5a} +  \nu_{5b}$ and $t_3 = 6 \nu_{5a} +  \nu_{5b} + 3
\nu_{7a} - 4 \nu_{7b}$. Since $|u^{7}|=5$, for any character
$\chi$ of $G$ we need to consider six cases, defined by part
(iv) of the Theorem.

Case 1. Let $\chi(u^{7}) = \chi(5a)$. By (\ref{E:2}) we obtain the system
\[
\mu_{5}(u,\chi_{2},*) = \textstyle \frac{1}{35} (4t_1 + 9) \geq
0; \qquad  \mu_{0}(u,\chi_{2},*)  = \textstyle \frac{1}{35} (-24t_1 + 16) \geq 0, \\ 
\]
which has no integral solutions.

Case 2. Let $\chi(u^{7}) = \chi(5b)$. Now the non-compatible inequalities are
\[
\mu_{7}(u,\chi_{2},*) = \textstyle \frac{1}{35} (6t_1 + 26) \geq
0; \qquad \mu_{0}(u,\chi_{2},*) = \textstyle \frac{1}{35} (- 24 t_1 + 36) \geq 0. \\ 
\]

Case 3. Let $\chi(u^{7}) = -2\chi(5a)+3\chi(5b)$.
By (\ref{E:2}) we obtain the system
\[
\mu_{7}(u,\chi_{2},*) = \textstyle \frac{1}{35} (6t_1 + 16) \geq
0; \qquad
\mu_{0}(u,\chi_{2},*)  = \textstyle \frac{1}{35} (-24t_1 + 76) \geq 0, \\ 
\]
which has no integral solution.

Case 4. Let $\chi(u^{7}) = -3\chi(5a)+4\chi(5b)$.
By (\ref{E:2}) we obtain the system
\[
\begin{split}
\mu_{0}(u,\chi_{2},*) & = \textstyle \frac{1}{35} (-24t_1 + 96)
\geq 0; \qquad \mu_{7}(u,\chi_{2},*) = \textstyle \frac{1}{35} (6t_1 + 11) \geq 0; \\ 
\mu_{0}(u,\chi_{3},*) & = \textstyle \frac{1}{35} (24t_2 + 175)
\geq 0; \qquad  \mu_{7}(u,\chi_{3},*)  = \textstyle \frac{1}{35}
(-6t_2 + 245) \geq 0; \\ 
\mu_{15}(u,\chi_{16},*) & = \textstyle \frac{1}{35} (4t_3 + 8071)
\geq 0; \qquad  \mu_{1}(u,\chi_{16},*) = \textstyle \frac{1}{35} (
-t_3 + 8001) \geq 0; \\ 
& \mu_{0}(u,\chi_{2},3)  = \textstyle \frac{1}{35} (-96 \nu_{5a} +
24 \nu_{5b} + 85) \geq 0, \\ 
\end{split}
\]
so $t_1=4$, $t_2 \in \{0, 35 \}$ and 
$t_3 \in \{ 21+35k \mid -58 \leq k \leq 228 \}$, and we have no solutions again.

Case 5. Let $\chi(u^{7}) = -4\chi(5a)+5\chi(5b)$. Using (\ref{E:2}), we obtain 
\[
\begin{split}
\mu_{7}(u,\chi_{2},*) & = \textstyle \frac{1}{35} (6t_1 + 6) \geq
0; \qquad \quad \; \mu_{1}(u,\chi_{2},*)  = \textstyle \frac{1}{35} (-t_1 - 1) \geq 0; \\ 
\mu_{0}(u,\chi_{3},*) & = \textstyle \frac{1}{35} (24t_2 + 155)
\geq 0; \qquad
\mu_{5}(u,\chi_{3},*)  = \textstyle \frac{1}{35} (-4t_2 + 155) \geq 0; \\ 
\mu_{15}(u,\chi_{16},*) & = \textstyle \frac{1}{35} (4t_3 + 8091)
\geq 0; \qquad
\mu_{1}(u,\chi_{16},*)  = \textstyle \frac{1}{35} ( -t_3 + 7996) \geq 0. \\ 
\end{split}
\]
Then $t_1=1$, $t_2 \in \{-5, 30 \}$ and
$t_3 \in \{ 16+35k \mid -58 \leq k \leq 228 \}$, and we have no solutions as
before.

Case 6. Let $\chi(u^{7}) = -\chi(5a)+2\chi(5b)$.
By (\ref{E:2}) we have incompatible inequalities
\[
\mu_{7}(u,\chi_{2}, *) = \textstyle \frac{1}{35} (6t_1 + 21) \geq
0; \qquad  \mu_{0}(u,\chi_{2},*)  = \textstyle \frac{1}{35}
(-24t_1 + 56) \geq 0. \\ 
\]

\noindent$\bullet$ Let $u$ be a unit of order $55$. By (\ref{E:1})
and Proposition \ref{P:4} we have
$$
\nu_{5a}+\nu_{5b}+\nu_{11a}+\nu_{11b}=1.
$$
Put $t_1 = 3 \nu_{5a} - 2 \nu_{5b}$, 
$t_2 = 6 \nu_{5a} + \nu_{5b}$ and  
$t_3 = 4 \nu_{5a} -  \nu_{5b} + 6 \nu_{11a} - 5 \nu_{11b}$.
Since  $|u^{11}|=5$, for any
character $\chi$ of $G$ we need to consider six cases, defined
by part (iv) of the Theorem.

Case 1. Let $\chi(u^{11}) = \chi(5a)$. 
Then by (\ref{E:2}) we obtain incompatible inequalities
\[
\mu_{5}(u,\chi_{2},*) = \textstyle \frac{1}{55} (4t_1 + 10) \geq
0; \qquad
\mu_{0}(u,\chi_{2},*) = \textstyle \frac{1}{55} (-40t_1 + 10) \geq 0. \\ 
\]

Case 2. Let $\chi(u^{11}) = \chi(5b)$. Using (\ref{E:2}) we obtain the system
\[
\begin{split}
\mu_{11}(u,\chi_{2},*) & = \textstyle \frac{1}{55} (10 t_1 + 20)
\geq 0; \qquad \; \mu_{0}(u,\chi_{2},*) = \textstyle \frac{1}{55}
(-40t_1 + 30) \geq 0; \\ 
\mu_{0}(u,\chi_{3},*) & = \textstyle \frac{1}{55} (40t_2 + 235)
\geq 0; \qquad  \mu_{11}(u,\chi_{3},*) = \textstyle \frac{1}{55}
(-10t_2 + 230) \geq 0; \\ 
\mu_{5}(u,\chi_{7},*) & = \textstyle \frac{1}{55} (4t_3 + 939)
\geq 0; \qquad \; \mu_{1}(u,\chi_{7},*) = \textstyle \frac{1}{55}
(-t_3 + 934) \geq 0, \\ 
\end{split}
\]
so
$t_1=-2$, $t_2\in \{ 1, 12, 23 \}$ and 
$t_3 \in \{ -1+55k \mid -4 \leq k \leq 17 \}$, and in every case 
we have no integer solution.

Case 3. Let $\chi(u^{11}) = -2\chi(5a)+3\chi(5b)$. By (\ref{E:2}) 
we obtain that
\[
\begin{split}
\mu_{11}(u,\chi_{2},*) & = \textstyle \frac{1}{55} (10t_1 + 10)
\geq 0; \qquad \; \mu_{0}(u,\chi_{2},*)  = \textstyle \frac{1}{55} (-40t_1 + 70) \geq 0; \\ 
\mu_{0}(u,\chi_{3},*) & = \textstyle \frac{1}{55} (40 t_2 + 195)
\geq 0; \qquad \mu_{11}(u,\chi_{3},*) = \textstyle \frac{1}{55}
(-10t_2 + 240) \geq 0; \\ 
\mu_{5}(u,\chi_{7},*) & = \textstyle \frac{1}{55} (4t_3 + 946)
\geq 0; \qquad \; \mu_{1}(u,\chi_{7},*)  = \textstyle \frac{1}{55}
(-t_3 + 891) \geq 0. \\ 
\end{split}
\]
From this follows that $t_1=-1 $, $t_2 \in \{2, 13, 24\}$ and 
$t_3 \in \{ 11+55k \mid -4 \leq k \leq 16 \}$, and for every
combination of $t_i$'s we have no solution.

Case 4. Let $\chi(u^{11}) = -3\chi(5a)+4\chi(5b)$. 
By (\ref{E:2}) we obtain the system
\[
\begin{split}
\mu_{11}(u,\chi_{2},*) & = \textstyle \frac{1}{55} (10 t_1 + 5)
\geq 0; \qquad  \mu_{0}(u,\chi_{2},*)  = \textstyle \frac{1}{55}
(-40 t_1 + 90) \geq 0, \\ 
\end{split}
\]
which has no integral solution.

Case 5. Let $\chi(u^{11}) = -4\chi(5a)+5\chi(5b)$.
By (\ref{E:2}) we have that
\[
\begin{split}
\mu_{11}(u,\chi_{2},*) & = \textstyle \frac{1}{55} (10 t_1 ) \geq
0; \qquad\qquad \; \mu_{1}(u,\chi_{2},*)  = \textstyle \frac{1}{55} (-t_1 ) \geq 0; \\ 
\mu_{0}(u,\chi_{3},*) & = \textstyle \frac{1}{55} (40t_2 + 155)
\geq 0; \qquad
\mu_{11}(u,\chi_{3},*) = \textstyle \frac{1}{55} (-10t_2 + 250) \geq 0; \\ 
\mu_{5}(u,\chi_{7},*) & = \textstyle \frac{1}{55} (4t_3 + 986)\geq
0; \qquad \; \mu_{1}(u,\chi_{7},*) = \textstyle \frac{1}{55}
(-t_3 + 881) \geq 0, \\ 
\end{split}
\]
so $t_1=0$, $t_2 \in\{3, 14, 25\}$ and 
$t_3 \in \{ 1+55k \mid -4 \leq k \leq 16 \}$, 
and again we have no solutions in every combination of $t_i$'s.

Case 6. Let $\chi(u^{11}) = -\chi(5a)+2\chi(5b)$. Using
(\ref{E:2}) we obtain the system
\[
\begin{split}
\mu_{11}(u,\chi_{2},*) & = \textstyle \frac{1}{55} (10t_1 + 15)
\geq 0; \qquad  \mu_{0}(u,\chi_{2},*)  = \textstyle \frac{1}{55}
(-40 t_1 + 50) \geq 0, \\ 
\end{split}
\]
which has no integral solution.


\noindent$\bullet$ Let $u$ be a unit of order $77$. 
By (\ref{E:1}) and Proposition \ref{P:4} we have
$$
\nu_{7a}+\nu_{7b}+\nu_{11a}+\nu_{11b}=1.
$$
Then using (\ref{E:2}) we obtain the non-compatible system of inequalities
\[
\begin{split}
\mu_{0}(u,\chi_{2},*) & = \textstyle \frac{1}{77} (60 (\nu_{7a} + \nu_{7b}) + 28) \geq 0; \\ 
\mu_{11}(u,\chi_{2},*) & = \textstyle \frac{1}{77} (-10(\nu_{7a} +\nu_{7b}) + 21) \geq 0. \\ 
\end{split}
\]
 
\bibliographystyle{plain}
\bibliography{Bovdi_Konovalov_McL}

\end{document}